\documentclass[12pt]{amsart}
\usepackage{amsmath,amsthm}

\def\pmatrix{\left(\begin{matrix}}
\def\endpmatrix{\end{matrix}\right)}

\def\P{{\mathbb P}}
\def\Z{{\mathbb Z}}
\def\C{{\mathbb C}}

\def\Q{{\mathbb Q}}

\def\cal{\mathcal}

\def\dd{{\rm d}}

\def\de{\delta}
\def\di{{\delta_{irr}}}
\def\slg2{\sum\limits_{i=1}^{\lfloor g/2\rfloor}}

\def\p{\partial}

\def\t{\theta}
\def\s{\sigma}
\def\T{\Theta}
\def\e{\varepsilon}

\def\w{\omega}

\def\A{{\mathcal A}}
\def\M{{\mathcal M}}
\def\J{{\mathcal J}}

\def\H{{\mathcal H}}

\def\ooo{\overline}

\def\omg{{\ooo{\M_g}}}
\def\mle{{\M_g^{2,4}}}
\def\ole{{\ooo\mle}}
\def\tt#1#2{{\t\left[\begin{matrix}{#1}\\{#2}\end{matrix}\right]}}

\def\WP{Weil-Petersson\ }

\theoremstyle{plain}
\newtheorem{thm}{Theorem}[section]
\newtheorem{lm}[thm]{Lemma}
\newtheorem{prop}[thm]{Proposition}
\newtheorem{cor}[thm]{Corollary}

\theoremstyle{definition}
\newtheorem{df}[thm]{Definition}
\newtheorem{rem}[thm]{Remark}
\newtheorem{ex}[thm]{Example}

\numberwithin{equation}{section}
\title{The degree of the Jacobian locus and the Schottky problem}
\author{Samuel Grushevsky}
\email{sam@math.princeton.edu}
\address{Mathematics Department, Princeton University,
Fine Hall, Washington Road, Princeton, NJ, 08540}
\date{December 28, 2003.}
\begin{document}

\begin{abstract}
We show that the degree of the images of the (level covers of)
moduli space of principally polarized abelian varieties $\A_g$ and
of the Jacobian locus $\J_g$ under the embedding $Th$ into the
projective space by theta constants are equal to the top
self-intersection numbers of one half the first Hodge class on
them. This allows us to obtain an explicit formula for $\deg
Th(\A_g)$ in all genera, compute $\deg Th(\J_g)$ for small $g$,
and obtain an explicit upper bound for $\deg Th(\J_g)$ for all
$g$.

Knowing $\deg Th(\A_g)$ allows us to effectively determine this
subvariety, i.e. effectively obtain {\it all} polynomial equations
satisfied by theta constants. Furthermore, combining the bound on
$\deg Th(\J_g)$ with effective Nullstellensatz allows us to
rewrite the Kadomtsev-Petvsiashvili (KP) partial differential
equation as a system of algebraic equations for theta constants,
and thus effectively obtain an algebraic solution to the Schottky
problem.
\end{abstract}
\maketitle

\section{Introduction}
The Schottky problem, the question of characterizing Jacobians of
Riemann surfaces among principally polarized abelian varieties
(ppavs), was solved by Shiota more than a hundred years after the
problem was first posed. Shiota completed the proof of Novikov's
conjecture: that a ppav is a Jacobian if and only if a certain
modification of its associated theta function satisfies the KP
equation \cite{shiota}. However, this solution of the Schottky
problem is not effective in the following sense: first, it
requires choosing the values of certain $3g+1$ parameters used to
modify the theta function, and, secondly, verifying that the KP
equation, a partial differential equation, holds for these values
of parameters. Checking the validity of such an equation is a hard
task --- it entails verifying that the sum of a certain Fourier
series is equal to zero at a given point, and this cannot be done
effectively.

Another approach to the Schottky problem, the Schottky--Jung
theory (see \cite{do},\cite{farkas}, and \cite{vangeemen} for a
review), if successful, would yield a system of algebraic
equations for theta constants of the second order defining the
locus of Jacobians inside the moduli space of ppavs, and thus
would give an effective algebraic solution to the Schottky
problem. However, so far this approach is only known to
characterize Jacobians up to additional components.

Numerous other approaches have been developed, but despite much
progress (see, for example, \cite{debarre} and Arbarello's
appendix to \cite{arbarello} for reviews), an entirely algebraic
solution to the Schottky problem, in the spirit of the Schottky's
original approach, has not yet been obtained for genus higher than
four. Some effective criteria to check that an abelian variety is
not a Jacobian have been obtained, though (see \cite{lazarsfeld}),
but if one produces an explicit ppav, there is no known way to
show that it actually is a Jacobian, apart from explicitly
constructing the corresponding curve.

In this paper we study the algebraic degrees of the image of the
Jacobian locus and of the moduli space of ppavs in the projective
space under the theta constants map, and obtain an effective
algebraic solution to the Schottky problem (``effective'' here
means that there is an actual explicit finite procedure that
produces the complete set of algebraic equations in theta
constants defining the Jacobian locus).

The structure of the work is as follows. We introduce notations in
section 2; in sections 3 and 4 we express the degrees in question
as some intersection numbers on moduli spaces. In section 5 we
compute and discuss these numbers for low genera. In section 6 we
use the methods of intersection theory on the moduli space,
together with Schumacher and Trapani's \cite{schum} and our
\cite{wpbound} results on Weil-Petersson volumes to obtain an
explicit upper bound for the degree of the Jacobian locus. In
sections 7 and 8 we then show how the KP equation can be
effectively rewritten as a system of algebraic equations for theta
constants characterizing Jacobians, and also how to effectively
obtain all relations among theta constants of general abelian
varieties.

\section{Notations}
Fix some genus $g>0$, which we will also interchangeably call the
dimension. Let $\Lambda$ be a lattice in $\C^g$ of maximal rank,
and let $A=\C^g/\Lambda$ be the associated abelian variety. We
will think of $A$ as being principally polarized and of $\Lambda$
as being generated over $\Z$ by the unit vectors in all directions
and the columns of a complex $g\times g$ matrix $\tau$, called the
period matrix for $A$. The period matrix $\tau$ has to be
symmetric with positive definite imaginary part; the set of such
$\tau$ is the {\it Siegel upper half-space} ${\cal H}_g$.
Different period matrices may define isomorphic ppavs --- this
corresponds to a change of basis in $\Lambda$ under the action of
the symplectic group $\Gamma_g:={\rm Sp}(2g,\Z)$ on the Siegel
upper half-space. Representing $\gamma\in\Gamma_g$ as
$\gamma=\pmatrix A&B\\ C&D\endpmatrix$, where $A,B,C,D\in{\rm
Mat}_{g\times g}(\Z)$, with the symplectic condition being
\begin{equation}
\pmatrix A&B\\ C&D\endpmatrix^t\pmatrix 0&1\\
-1&0\endpmatrix \pmatrix A&B\\ C&D\endpmatrix=\pmatrix 0&1\\
-1&0\endpmatrix,
\end{equation}
the action of $\Gamma_g$ on ${\cal H}_g$ is given by
$\gamma\tau:=(A\tau+B)(C\tau+D)^{-1}$. We denote by $\A_g$ the
fundamental domain for this action, which is the moduli space of
ppavs of dimension $g$ up to biholomorphisms: $\A_g:={\cal
H}_g/\Gamma_g$.

We also define {\it level} $n$ subgroups of the symplectic group,
for $n$ a positive integer:
\begin{equation}
\begin{matrix}
\Gamma_g(n):=\left\lbrace \gamma=\pmatrix A&B\\
C&D\endpmatrix \right|\left. \gamma\in\Gamma_g, \pmatrix A&B\\
C&D \endpmatrix \equiv \pmatrix 1&0\\ 0& 1\endpmatrix {\rm\
mod\ } n\right\rbrace\\
\Gamma_g(n,2n):=\left\lbrace \gamma \right|\left. \gamma\in
\Gamma_g(n), {\rm diag}(AB^t)\equiv{\rm diag}(C D^t)\equiv 0
{\rm\ mod\ } 2n \right\rbrace.
\end{matrix}
\end{equation}

A function $f:{\cal H}_g\to \C$ is called a {\it modular form of
weight $k$ with respect to a subgroup $\Gamma\subset\Gamma_g$} if
\begin{equation}
f(\gamma\tau)=\det (C\tau +D)^kf(\tau)\qquad
\forall\gamma\in\Gamma ,\ \forall\tau\in{\cal H}_g.
\end{equation}

We define {\it theta functions of the second order} to be, for
$\e\in \frac{1}{2}\Z^g/\Z^g$, $\tau\in{\cal H}_g$, and $z\in\C^g$
\begin{equation}
\Theta[\e](\tau,z):=\sum\limits_{n\in\Z^g}\exp(2\pi
i(n+\e)^t\tau(n+\e) +2\pi i (n+\e)^tz).
\end{equation}
The values of theta functions for $z=0$ are called {\it theta
constants}. These theta constants are modular forms of weight one
half with respect to  $\Gamma_g(2,4)$ (see \cite{igusa}).
Geometrically they are sections of a certain line bundle $L$ over
a finite cover $\A_g^{2,4}:={\cal H}_g/ \Gamma_g(2,4)$ of $\A_g$.
Under the action of $\Gamma_g(2,4)$ all theta constants multiply
by the same factor $\det (C\tau +D)^{1/2}$, and thus the mapping
\begin{equation}
Th:\tau\mapsto\lbrace\Theta[\e](\tau) \rbrace_{\e\in\frac{1}{2}\Z^g/\Z^g}
\end{equation}
gives a well-defined map $Th:\A_g^{2,4}\to\P^{2^g-1}$, which is
known to be generically injective (see \cite{sasaki},
\cite{smlevel}). In \cite{smlevel} Salvati Manni claims the
injectivity of the map, but he has informed us that there is a
small gap in the proof. For our purposes here generic injectivity
combined with the statement that $Th$ is at most finite-to-one
everywhere, (which follows from the fact that the theta-null map
is an embedding, see \cite{igring}) suffices. Since
$\Gamma_g(2,4)$ is a normal subgroup of $\Gamma_g$, the level
moduli space $\A_g^{2,4}$ is a Galois cover of $\A_g$, of degree
(see \cite{igring}) equal to
\begin{equation}
[\Gamma_g:\Gamma_g(2,4)]=2^{g^2+2g}\prod\limits_{k=1}^g
(2^{2k}-1). \label{degcover}
\end{equation}

We denote by $\M_g$ the {\it moduli of Riemann surfaces of genus
$g$}, by $J:\M_g\to\A_g$ --- the {\it Torelli map} sending a curve
$X$ to its Jacobian $J(X)$, and by $\J_g:=J(\M_g) \subset\A_g$ the
locus of (isomorphism classes of) Jacobians within the moduli of
ppavs. Let $\A_g^{\rm irr}$ be the locus of irreducible ppavs;
then in fact the image $J(\M_g)\subset\A_g^{\rm irr}$.

Corresponding to the level cover $p:\A_g^{2,4}\to\A_g$ we can
consider the level cover $p:\M_g^{2,4}\to\M_g$ with the level
Torelli embedding $J^{2,4}:\M_g^{2,4}\to\A_g^{2,4\rm irr}$, whose
image we denote by $\J_g^{2,4}$ (see \cite{vangeemen},\cite{do}).
Schottky problem is the question of describing
$Th(\J_g^{2,4})\subset Th(\A_g^{2,4})\subset\P^{2^g-1}$.

Recall that the moduli spaces $\M_g$ and $\A_g$ are not compact,
and this is where most of our trouble will come from. We will need
to use the Deligne-Mumford compactification (see
\cite{demu},\cite{harris}) $\omg$ of the moduli space of curves,
and a certain second Voronoi toroidal compactification
$\ooo{\A_g}$ of $\A_g$, such that the embedding $J:\M_g\to \A_g$
extends to a map $\ooo J:\omg\to\ooo{\A_g}$, injective on the
locus of irreducible stable curves, and which becomes an embedding
if we blow up the locus of reducible abelian varieties in $\A_g$
(see \cite{namikawa}, chapter 9.D). There also exists another
compactification $\widetilde\A_g$ of $\A_g$, called the Satake, or
minimal, compactification, from which the toroidal
compactification can be obtained by a series of blowups. In
particular there is a natural blow-down map
$b:\ooo\A_g\to\widetilde\A_g$, and the theta constants extend
naturally to define a map $Th:\widetilde\A_g\to\P^{2^g-1}$.
Abusively we denote the compose map $b\circ Th:\ooo{\A_g^{2,4}}\to
\P^{2^g-1}$ also by $Th$. We let $\ole$ be the Deligne-Mumford
compactification of the level moduli space (i.e. the normalization
of $\omg$ in the field of functions of $\mle$) --- see
\cite{demu}, \cite{bernstein} for more details. We also have the
extended level Torelli map
$\ooo{J^{2,4}}:\ole\to\ooo{\A_g^{2,4}}$.

\section{Theta metric}
\begin{df} We call the {\it theta metric} the restriction of
the Fubini-Study metric from the cotangent bundle of $\P^{2^g-1}$
to the cotangent bundle of the image $Th(\A_g^{2,4})$. When we
pull it back to $\A_g^{2,4}$, it is a metric on the bundle $L$ of
theta constants. By abuse of notation we will also call theta
metric the pullback of this metric to $\M_g^{2,4}$ by the map
$J^{2,4}$.
\end{df}

Denote by $\Omega$ the curvature form of the Fubini-Study metric
on $\P^{2^g-1}$ and, by abuse of notation, its pullbacks to
$\A_g^{2,4}$ and $\mle$. It is given by (see \cite{GrHa}, section
0.2)
\begin{equation}
\Omega=\frac{\sqrt{-1}}{2\pi}\p\ooo\p\log
\sum\limits_{\e\in\Z_2^g}|\T[\e](\tau)|^2.
\end{equation}

By definition the degree of
$\overline{Th(\J_g^{2,4})}\subset\P^{2^g-1}$ is equal to the
integral of the top power of the form $\Omega$ over it:
\begin{equation}
\deg \overline{Th(\J_g^{2,4})}=
\int\limits_{Th(\J_g^{2,4})}\Omega^{3g-3}.
\end{equation}
Note that when computing the integral, we can as well integrate
over the open part, while technically the degree only makes sense
for closed subvarieties. As we do not know the image
$Th(\J_g^{2,4})$ explicitly (knowing it would be solving the
Schottky problem), explicitly integrating over it is unmanageable.
Thus we pull $\Omega$ back to $\mle$, and try to integrate there.
However, integrating over $\mle$ is still unmanageable, and
instead we express the integral as an intersection number on the
compactification $\ole$. But we do not know the intersection
theory on $\ole$ well enough to be able to compute this
intersection number. Using the branched cover $\ooo
p:\ole\to\omg$, we reduce the intersection-theoretic computation
on $\ole$ to an intersection-theoretic computation on $\omg$,
which is then manageable. In doing all of this, the main
difficulty is dealing with the boundary contributions
appropriately.

\begin{rem}
It is known that, up to torsion, over $\A_g^{2,4}$ the bundle $L$
of modular forms of weight one half is the square root of the
determinantal line bundle of the Hodge bundle $H$ --- the
$g$-dimensional complex vector bundle with the fiber over some $A$
being $H^{1,0}(A,\C)$, the space of holomorphic one-forms on $A$
(see \cite{faltings}, section 1.5). We would then like to say that
the integral of the top power of $\Omega$ is just equal to the top
self-intersection number of the divisor corresponding to the
determinantal line bundle of the Hodge vector bundle. This indeed
happens to be the case, but actually {\it proving} it requires
some work, as the theta metric has a singularity near the boundary
of the and it is not clear a priori how the bundle and the metric
extend to the boundary. One possible way to deal with the problem
is to extend the bundle and then smoothen the metric using
complex-analytic techniques of Monge-Amp\`ere operator applied to
currents (see \cite{bt}) to compute its total mass (see
\cite{ra}). What we will do for our proof, however, is study the
boundary degeneration directly. In \cite{geom} Wolpert deals with
the similar problem for the \WP metric.
\end{rem}

\begin{rem}
There is a natural metric on the Hodge bundle over $\M_g$, called
the Hodge metric, discussed by Nag in \cite{nag}: for two abelian
differentials $\w$ and $\w'$ on $X$ we define their scalar product
to be $\int_X\w\wedge\ooo{\w'}$. We will show that Hodge and theta
metrics are distinct by showing that their behaviors on the
boundary are different. We will obtain explicit upper bounds for
theta volumes, while at this time we were not able to obtain a
bound for Hodge volumes or to express them as intersection
numbers.
\end{rem}

\begin{ex}
To illustrate the potential problems arising, we provide a simple
example (due to Yum-Tong Siu) illustrating the difference between
the integral of the top power of a current and the top
self-intersection number of the cohomology class it represents.

Indeed, on $\P^2$ with homogeneous coordinates $(x:y:z)$ consider
the expression $\p\ooo\p\log(|x|^2+|y|^2)$ --- it is a smooth form
outside the point $(1:0:0)$ and can be trivially extended to a
current $\phi$ on $\P^2$. The class $[\phi]\in H^2(\P^2)=\Q$ is
clearly non-zero (in fact it is equal to the class of the
Fubini-Study metric), and thus the intersection number $\langle
[\phi]^2\rangle_{\P^2}\ne 0$. On the other hand the integral
$\int_{\P^2\setminus(1:0:0)} \phi^{\wedge 2}$ is zero, since the
form $\phi$ is exact on the chart where $x\ne0$.
\end{ex}

\begin{prop}
The form $\Omega$ extends smoothly to the compactifications $\ole$
and $\ooo{\A_g^{2,4}}$. We denote these extensions by
$\ooo\Omega$.
\end{prop}
\begin{proof}
Theta constants extend holomorphically to the Satake
compactification $\widetilde{\A_g^{2,4}}$ (see \cite{faltings},
chapter V). Thus the curvature form $\Omega$ extends smoothly to
$\widetilde\Omega$ on the Satake compactification. Then
$\ooo\Omega:=b^*(\widetilde\Omega)$ is a smooth extension of
$\Omega$ to $\ooo{\A_g^{2,4}}$.
\end{proof}

\section{Theta constants at the boundary}
\begin{df} We denote by $D_i$ the boundary divisor of $\omg$
consisting of the closure in $\omg$ of the locus of stable curves
that are bouquets of two curves of genera $i$ and $g-i$
respectively. We also denote by $D_0$ the (closure of) the locus
of irreducible nodal curves in $\omg$. By Poincar\'e duality on
$\omg$ each of these divisors determines a class in $H^2(\omg)$.
We denote the corresponding cohomology classes by $\de_i$ and
$\di$, respectively. Let us also denote by $\lambda$ the {\it
first Hodge class}, i.e. the first Chern class of the Hodge bundle
on $\omg$.
\end{df}

We would now want to say that the divisor of the bundle $L$ of
modular forms of weight one half is equal to the divisor of some
specific theta constant, which is its section, and thus that the
divisor of one theta constant is equal to the divisor of curves
which have {\it some} theta constant vanishing, divided by $2^g$.
However, as pointed out to us by Carel Faber, this argument does
not necessarily have to work, and this is where we have to worry
about level structures. The problem is that the divisor of one
theta constant is not defined on $\omg$, and that the covering map
$\ooo p$ is branched exactly over the boundary and over the locus
of reducible abelian varieties: see \cite{do},\cite{vangeemen} and
\cite{smlevel}. However, this argument can be salvaged and leads
to the solution.

\begin{prop} The form $\ooo\Omega$ is invariant under the change of
level structure, i.e. under the deck transformations of the cover
$\ooo p:\ooo{\A_g^{2,4}}\to\ooo{\A_g}$. Since the cover $\ooo p$
is unbranched over $\A_g^{\rm irr}$, see \cite{smlevel}, on
$\A_g^{2,4\rm irr}$ the form $\Omega$ is a pullback of a smooth
form on $\A_g^{\rm irr}$, $\Omega=p^*\omega$ for some $\omega\in
H^{1,1}(\A_g^{\rm irr})$. This also holds for $\M_g$ and we
abusively use the same notations there.
\end{prop}
\begin{proof}
The proof easily follows from the transformation laws for theta
functions (see, for example \cite{bila}, 8.6.1). A very detailed
study of the transformations of theta constants is done in
\cite{smprojvar}, \cite{smlevel}, and the following proof is
basically a quote from there.

The action of the deck transformations group
$\Gamma_g/\Gamma_g(2,4)$ on the set of theta constants of the
second order is as follows:

--- if $\gamma\in\Gamma_g(2)/\Gamma_g(2,4)$, then (apart from the modular
factor that is the same for all theta constants, and thus would
not matter for the Fubini-Study metric) $\gamma$ acts on each
$\T[\e]$ by multiplying it by a sign, which may be different for
different $\e$, but still does not matter for computing $\Omega$.

--- if $\gamma\not\in\Gamma_g(2)/\Gamma_g(2,4)$, then in addition to all of
the above, the action of $\gamma$ permutes theta constants
$\T[\e]$ in a certain way, but this again does not change the form
$\Omega$.

Thus the form $\ooo\Omega$ on $\ooo{\A_g^{2,4}}$ is invariant
under the deck transformations. Thus over the locus where the
cover is unbranched it is the  pullback of some smooth form, i.e.
on $\A_g^{2,4\rm irr}$ it must be the (invariant) pullback of some
form $\omega$ on $\A_g^{\rm irr}$.
\end{proof}

\begin{prop} Extend the form $\omega$ trivially to a current on
$\ooo{\A_g}$, and denote by $[\omega]\in H^{1,1}(\ooo{\A_g})$ its
cohomology class. The class $[\ooo\Omega]\in
H^{1,1}(\ooo{\A_g^{2,4}})$ is invariant under the deck
transformations, and $[\ooo\Omega]=p^*[\ooo\omega]$.
\end{prop}
\begin{proof}
Since the form $\Omega$ is invariant under deck transformations,
the divisor class $[\Omega]$ it represents is also invariant, and
is the pullback of the divisor class $[\omega]$ of the form of
which $\Omega$ is the pullback. This works perfectly when the
cover is not branched. Over the branching locus we still have the
invariance of $\ooo\Omega$ under the deck transformations, but now
$\ooo\omega$ may actually have singularities, while its pullback
$\ooo\Omega$ stays smooth. However, the above statement about
divisor classes still makes sense: the closed currents also
represent cohomology classes, and we just need to compute these.
Finally this means that we still have
$[\ooo\Omega]=p^*[\ooo\omega]$.
\end{proof}

\begin{rem}
The local picture near the branching locus is as follows: suppose
locally the branching order is $B$, so that in local coordinates
the covering map $p:z\to w$ is $w=p(z)=z^B$. Then for any $\e\ge0$
the pullback of the current $|w|^{-2+2/B+\e}\dd w\wedge\dd\ooo w$
is the smooth form
$$
B^2|z|^{-2B+2+B\e}z^{B-1}\ooo{z}^{B-1} \dd z \wedge\dd\ooo{z}=
B^2|z|^{B\e}\dd z\wedge\dd\ooo z.
$$
\end{rem}

Combining all the equalities between integrals and intersection
numbers that we have obtained and using the fact that $\Omega$ is
smooth on $\ooo{\A_g^{2,4}}$, we get
\begin{thm}
The degree of the Jacobian locus is
\begin{equation}
\deg\ooo{Th(\J_g^{2,4})}=\int\limits_\mle \Omega^{3g-3}=
\left\langle[\ooo\Omega]^{3g-3}\right\rangle_{\ole}= \deg
p\,\left\langle[\ooo\omega]^{3g-3}\right\rangle_\omg
\end{equation}
and similarly for the degree of $\ooo{Th(\A_g^{2,4})}$.
\label{fformula}
\end{thm}
Thus finally we are left with the question of computing the class
$[\ooo\omega]\in H^{1,1}(\omg)$ and also of $[\ooo\omega]\in
H^{1,1}(\ooo{\A_g})$.

\smallskip
Before proceeding to do this in the next section, let us give an
analytic description of the picture for the moduli of curves. The
complex coordinates near the boundary of the moduli space are
given by the plumbing construction (introduced by Bers in
\cite{bers}, see also \cite{masur},\cite{geom}): we take a stable
curve $X$, cut out a small neighborhood of every node $p_i$, which
locally looks like $\lbrace zw=0\rbrace$ in $\C^2$, and replace it
by the neighborhood $\lbrace zw=s_i\rbrace$ for some $s_i\in\C$.
Then $\lbrace s_i\rbrace $ and the rest of the complex coordinates
$\lbrace t_i\rbrace$ on the moduli given by the Bers embedding and
deformation theory give complex coordinates in a neighborhood of
this noded curve.

In these coordinates, the asymptotic behavior of the period matrix
was obtained by Fay and by Yamada (see \cite{yamada},
\cite{taniguchi}). In particular, if $X$ is an irreducible stable
curve with one node, i.e. $X\in D_0$, then in a neighborhood of
$X$ the period matrix stays bounded except for $\tau_{11}$, which
grows as $\frac{\log s_1}{2\pi i}$ as $s_1\to 0$.

Using the Fourier-Jacobi expansion (see \cite{faltings}, chapter
V) for theta functions
$$
\T[\e_1\ \e']\left(\begin{matrix}\tau_{11} &\xi^t\\
\xi &\tau'\end{matrix}\right)=\sum\limits_{m\in\Z}\exp\left(4\pi i
\tau_{11} (2m+\e_1)^2\right) \T[\e'](\tau',(m+\frac{\e_1}{2})\xi),
$$
we can then compute the asymptotics of the growth of the form
$\omega$ near the boundary $D_0$ of $\omg$ --- this computation is
very similar to the one in \cite{do}, and we omit the easy
details. The result is given by
\begin{prop}
The form $\omega$ has a singularity of type $|s_1|^{- 3/2}\dd
s_1\wedge \dd\ooo{s_1}$ near the boundary component
$D_0\subset\omg$. The volume form $\omega^{\wedge(3g-3)}$ has a
singularity of the same order.
\end{prop}

\begin{rem}
For comparison, for the Weil-Petersson volume form the singularity
is (see \cite{masur}) $|s_1|^{-2}(\log|s_1|)^{-3}\dd
s_1\dd\ooo{s_1}$, and for the curvature of the Hodge metric the
singularity is $|s_1|^{-2}(\log|s_1|)^{-2}\dd s_1\dd\ooo{s_1}$.
The last fact is easy and well-known; to see it we just notice
that the Hodge metric is $\det\operatorname{Im}\tau\sim \log
|s_1|$, and thus its curvature form blows up as
$$
\p\ooo\p\log \det\operatorname{Im}\tau \sim\p\ooo\p\log(\log
|s_1|).
$$
Thus the Hodge curvature has the worst singularity, while $\omega$
has the mildest singularity of the three.

Notice also that for the Weil-Petersson and Hodge metrics
$|s_1|^{-1}$ is taken in the power 2, while for $\omega$ it is
taken to the power 3/2. This indicates that while potentially
there can be delta-function contributions of the boundary to the
volume computation for Hodge and Weil-Petersson metrics, for the
theta metric the boundary can be excluded, as the singularity
there is only branching, as we already know. As we will explain
below, this basically translates to the fact that the divisor
class of $[\omega]$ lives in $H^2(\M_g)$, {\it uncompactified},
i.e. does not have any extra terms coming from the boundary
divisors.
\end{rem}

\section{The degree computations}
\begin{thm}
In $H^2(\omg)$ we have $[\ooo\omega]=\lambda/2$. Combined with
theorem \ref{fformula}, this means that the degree of the Jacobian
locus is
$$
\deg\ooo{Th(\J_g^{2,4})}=\deg p\,\left\langle(\lambda/2)^{3g-3}
\right\rangle_{\omg}.
$$
\end{thm}
\begin{proof}
Analytically this follows from the fact that the singularity of
the current $\ooo\omega$ near the boundary of $\omg$ is milder
than $|z|^{-2}\dd z\wedge\dd \ooo z$, so that there is no delta
function on the boundary (i.e. in some sense the integral of the
volume form $\ooo\omega^{3g-3}$ over the boundary is zero) and
thus there are no $\delta_{\rm irr}$ summands in $[\ooo\omega]$.
The fact that there are no $\delta_i$ summands in $[\ooo\omega]$,
either, follows similarly from the smoothness of the form
$\ooo\Omega$ on $\A_g^{2,4}$
--- note that similarly to the behavior near $\p\ooo{\A_g}$, the
form $\ooo\omega$ may blow up near $D_i$, as the cover $\ooo p$ is
branched there. To make this proof rigorous, however, one needs to
also study the vanishing near higher-order degenerations and,
though doable, this becomes more cumbersome, unless we try to
quote directly the fact that any branching-type singularities do
not matter for computing the total mass of a current, for which we
were not able to find a reference.

Thus let us present a rigorous proof specific to the problem at
hand. Basically it stems from the fact that the form $\Omega$
extends smoothly to the Satake compactification of $\M_g$ (or of
$\A_g$), and there the boundary is not a collection of divisors,
but higher codimension, and thus there are no ``boundary classes''
in $H^2$. Indeed, let us use the basis $\lbrace{\cal E}, {\cal
L}_i\rbrace$ for $H_2(\omg)$ constructed by Wolpert in
\cite{wolpert}: here ${\cal E}$ is the family of elliptic tails,
i.e. a varying 1-pointed elliptic curve attached to a fixed genus
$g-1$ curve with 1 marked point, and ${\cal L}_i$ are obtained by
attaching a varying 4-pointed sphere to one or two fixed Riemann
surfaces of lower genera.

Then by definition we have
$$
[\ooo\omega]([{\cal L}_i])_\omg= [\ooo\omega]([\ooo J({\cal
L}_i)])_{\ooo{\A_g}}=(\deg p)^{-1}[\ooo\Omega]([p^*\circ\ooo
J({\cal L}_i)])_{\ooo{\A_g^{2,4}}}
$$
$$
=(\deg p)^{-1}[\widetilde\Omega]\left([b\circ p^*\circ \widetilde
J({\cal L}_i)]\right)_{ \widetilde{\A_g^{2,4}}}=(\deg
p)^{-1}[\widetilde\Omega]([0])=0
$$
because the Jacobian of a bouquet of smooth curves does not keep
track of the point of attachment, and, more generally, the image
of the level cover of the family ${\cal L}_i$ (in which the curves
differ only by the choice of points where the fixed components are
attached to the sphere), is locally constant in the level cover of
Satake compactification, i.e. is simply a collection of points.

Checking the intersection matrix of $H^2(\omg)$ with $H_2(\omg)$
in \cite{wolpert}, we see that it follows that $[\ooo\omega]$ is
proportional to $\lambda$. Then either computing the intersection
with $\cal E$ or recalling that over $\M_g$ the theta constant
bundle is the square root of the determinant of Hodge bundle we
see that in fact $[\ooo\omega]=\lambda/2\in H^2(\omg)$.
\end{proof}

The top self-intersection numbers of the class $\lambda$ on $\omg$
were computed by Carel Faber using his program described in
\cite{fabnum}. For the degrees of the Jacobian locus in genera
$1\ldots 7$ one then gets
\begin{equation}
\begin{matrix}1,1,16,208896,282654670848,\\
23303354757572198400,87534047502300588892024209408.
\end{matrix}
\end{equation}
These numbers were known before in genera one through three: in
genera one and two there are no equations by dimension
considerations, and the answer of 16 (obtained by rewriting the
Schottky's original equation in terms of $\T[\e]$) for genus three
is explained in \cite{vgvdg}. The number in genus four is already
large and thus rather worrying. However, this is the ``total
degree'' of all equations needed to define
$Th(\J_4^{2,4})\subset\P^{15}$, and not only the degree of the
additional equations defining $Th(\J_4^{2,4})\subset
Th(\A_4^{2,4})$.

\begin{thm}
The later degree, $\deg Th(\ooo{\A_g^{2,4}})
=\deg\ooo{Th(\A_g^{2,4})}$, is equal to $\deg p$ times the top
self-intersection number of $\lambda/2$ on $\ooo\A_g$.
\end{thm}
\begin{proof}
Analytically this works the same way as for the moduli of curves
in the previous theorem, so this is again not very rigorous. A
rigorous proof of the fact that the top self-intersection of
$\lambda$ does not include anything on $\p\ooo{\A_g}$ is
essentially contained in \cite{vdgeer}.

A self-contained rigorous proof is as follows: it is known from
\cite{mumford} (see also \cite{hulek} for a review and an
interesting discussion of related issues) that the Picard group of
$H^2(\ooo\A_g)$ is generated by $\lambda$ and the class of the
boundary $\delta$. To show that the class $[\ooo\omega]$ is
proportional to $\lambda$, i.e. does not contain a boundary
summand, we can use the basis $\ooo J(\cal E)$, $\ooo J({\cal
L}_0)$ for $H_2(\ooo\A_g)$ --- the fact that this is a basis
follows from the non-degeneracy of the $2\times 2$ intersection
matrix of these two homology classes with $\lambda,\delta$, a
basis for $H^2(\ooo\A_g)$ --- and note that the intersection of
$[\ooo\Omega]$ with $\ooo J({\cal L}_0)$ is zero, as above.
\end{proof}

\begin{rem} Corrado De Concini explained to us that the last
theorem can also be proven algebraically in a different way.
Indeed, the degree of a variety is the number of points in its
intersection with a generic linear subspace of complementary
dimension. Since the boundary of the Satake compactification has
high codimension, such a linear subspace can be chosen not to to
intersect the boundary, so the boundary would not matter, and we
can compute the intersection number on the Satake compactification
(where there is only one divisor), which it follows from the work
\cite{vdgeer} is the same. Such an argument does not work for the
moduli space of curves, as there we do not know how to compute
intersection numbers on the Satake compactification.
\end{rem}

Carel Faber has brought the work of van der Geer \cite{vdgeer} to
our attention. One of the results of that paper is a closed
formula for the top self-intersection of $\lambda$ on $\A_g$:
\begin{prop}[\cite{vdgeer}]
Denoting $N:=g(g+1)/2$, in the orbifold sense (which, as explained
to us by Carel Faber, means we will later need to multiply these
by two, due to the presence of $x\to-x$ involution on every
abelian variety) the intersection numbers we get are
\begin{equation}
\langle \lambda^N\rangle_{\ooo\A_g}=\langle
\lambda^N\rangle_{\A_g}=(-1)^N
N!\prod\limits_{k=1}^g\frac{\zeta(1-2k)}{2\left((2k-1)!!\right)},
\end{equation}
where $\zeta(1-2k)$ are the values of Riemann's zeta function at
negative odd integer points (equal to some products of Bernoulli
numbers and factorials), and double factorial denotes the product
$(2k-1)!!:=(2k-1)(2k-3)(2k-5)\cdots 3\cdot 1$.
\end{prop}

Multiplying these intersection numbers by $2(\deg p)/2^N$ we can
compute $\deg Th(\ooo{\A_g^{2,4}})$ for any $g$. In particular for
genera $1\ldots 7$ we get
\begin{equation}
\begin{matrix}
1,1,16,13056,1234714624,\\
25653961176383488,197972857997555419746140160.
\end{matrix}
\end{equation}
In genera one and two this agrees with the fact that there
are no equations, and the only relation in genus 3 has degree
16 (see \cite{vgvdg}, p.~623).

We also compute the ratio $\frac{\deg Th(\J_g^{2,4})}{\deg
Th(\A_g^{2,4})}$ for genera $1\ldots 7$:
\begin{equation}
1,1,1,16,\frac{2976}{13},\frac{202742400}{223193},
\frac{8678490624}{19627855}.
\end{equation}
The 1's we get for genera one through three correspond to the fact
that $\M_g=\A_g^{\rm irr}$ for $g\le 3$. For genus four $Th(\J_4)$
is defined within $Th(\A_4)$ by one extra equation, Schottky's
original relation, which is of degree 16 (see \cite{vgvdg}), so
our computation produces the correct result. The fact that the
ratios we get in higher genera are not integer yields the
following
\begin{cor}
The variety $\ooo{Th(\J_g^{2,4})}$ is not a complete intersection
within $\ooo{Th(\A_g^{2,4})}$ for genera 5,6,7.
\end{cor}

\begin{rem}
There is a nuance here: in genus four the defining equation is of
degree 16 in theta-nulls, not in theta constants of the second
order. To rewrite it algebraically in terms of $\T$'s, it has to
be multiplied by three conjugates and becomes of degree 64 (see
\cite{vgvdg}), but only one of the four isomorphic components of
the locus it then defines is the Schottky locus. It is also
possible to perform computations analogous to ours for the
embedding of $\A_g^{4,8}$ by theta-nulls --- we will then get a
larger factor corresponding to the degree of the level (4,8) cover
instead of the level (2,4) cover, but the ratio of the degrees of
the Jacobian locus and of the moduli of abelian varieties will
stay the same, 16. We will study $Th(\A_4^{2,4})$ and
$Th(\J_4^{2,4})$ in more detail in another work.
\end{rem}

Our degree computations also allow us to prove one curious
result:
\begin{prop}
Consider the (level cover of) the locus of 1-reducible abelian
varieties: $\A_1^{2,4}\times\A_{g-1}^{2,4}\subset\A_g^{2,4}$. Its
image under the theta map is not a complete intersection in
$Th(\A_g^{2,4})$ for $g=4,5,6,7$. Similarly for the moduli space
of curves $Th(J^{2,4}(D_1)))$ is not a complete intersection in
$Th(\J_g^{2,4})$ for $g=4,5,6,7$.
\end{prop}
\begin{proof}
To prove this we note that from the considerations above it
follows that
$$
\deg Th(\A_1^{2,4}\times\A_{g-1}^{2,4})=\deg
p\left\langle\frac{\lambda}{2}\right\rangle_{
\ooo{A_1}}\left\langle\left(\frac{\lambda}{2}\right)^{g(g-1)/2}
\right\rangle_{\ooo{A_{g-1}}},
$$
and, evaluating the intersection number on $\A_1$, that
$$
\frac{\deg Th(\A_1^{2,4}\times\A_{g-1}^{2,4})}{\deg
Th(\A_g^{2,4})}=\frac{2^{g-1}\left\langle\lambda^{g(g-1)/2}
\right\rangle_{\ooo{A_{g-1}}}}{24\left\langle\lambda^{g(g+1)/2}
\right\rangle_{\ooo{A_g}}}
$$
and similarly for the moduli spaces of curves. These ratios are
not integer for the values of $g$ in question, and thus one
variety cannot be a complete intersection inside the other.
\end{proof}

This proposition is of interest for the following reason: if the
locus of reducible abelian varieties were a complete intersection,
i.e. were given inside $Th(\A_g^{2,4})$ by $g-1$ equations, then
restricting this to the Jacobian locus we could potentially get a
cover of $\M_g^{2,4}$ by $g-1$ affines, each corresponding to the
non-vanishing of one of the equations, and would prove Looijenga's
conjecture that the cohomology of the uncompactified $\M_g$ has
properties similar to that of a $(g-1)$-dimensional variety. Thus
the proposition shows that this approach to proving Looijenga's
conjecture is basically hopeless.

\section{Intersection number estimates}
We have reduced the computation of the degree of the locus of
Jacobians to the computation of the theta volumes, which we then
expressed as an intersection number on $\omg$. Now we will obtain
an explicit upper bound for these intersection numbers.

\begin{thm}
The degree of the Jacobian locus can be bounded explicitly:
$$
\deg Th(\ooo{\J_g^{2,4}})<C(g),
$$
where $C(g)$ is an explicit function of $g$ with the leading
growth order $2^{2g^2}$. \label{dbound}
\end{thm}
To obtain the bound we will use the results from previous work on
Weil-Petersson volumes. Let $\kappa$ be the first Mumford's
tautological cohomology class on $\M_{g,n}$. For us it is just the
class of the K\"ahler form of the \WP metric, divided by $2\pi^2$
(see \cite{wolpert}): $\kappa=[\omega_{WP}]/2\pi^2$. Then
algebraically instead of looking at ${\rm Vol}_{g,n}$, the
Weil-Petersson volume of $\M_{g,n}$ (see \cite{wpbound}), we can
look at the top self-intersection number
\begin{equation}
\langle\kappa^{3g-3+n}\rangle_{\ooo{\M_{g,n}}}=\frac{(3g-3+n)!{\rm
Vol}_{g,n}} {(2\pi^2)^{3g-3+n}}.
\end{equation}
\begin{rem}
This equality is not at all trivial: the \WP form $\omega_{WP}$ is
only smooth on $\M_g$, and extends to $\omg$ as a closed positive
current, and we get the same kind of problem as we had when
extending $\Omega$. The \WP volume is the total Monge-Amp\`ere
mass (see \cite{bt}, \cite{ra}) of this current, which in general
might not be equal to the top self-intersection number of its
cohomology class. The equality, proven in \cite{geom} and further
discussed in \cite{bundle}, follows from the fact that the form
$\omega$ extends to the boundary smoothly in the smooth structure
near the boundary defined by certain Fenchel-Nielsen coordinates,
and the extension has the same \v Cech cohomology class as the
trivial extension of $\omega$ as a current in complex coordinates.
\end{rem}

In \cite{wpbound} we used Penner's decorated Teichm\"uller theory
and his earlier work \cite{wpvol} on the subject to show that for
$n$ fixed and $g$ large ${\rm Vol}_{g,n}<c^g g^{2g}$ for some
explicit constant $c$. In \cite{schum} Schumacher and Trapani use
ampleness of $\kappa$ (which follows from the \WP metric being
K\"ahler and extendable smoothly to $\omg$), and the knowledge of
some explicit effective divisors on $\ooo{\M_{g,n}}$ to obtain
lower bounds on \WP volumes. Indeed, since $\kappa$ is ample on
$\omg$, we have $\langle\kappa^{3g-4}D\rangle_\omg\ge 0$ for any
effective divisor $D$. Since the \WP metric on the boundary of the
moduli restricts to the \WP metric on the lower-dimensional moduli
spaces, choosing some particular $D$'s then yields estimates on
${\rm Vol}_{g,n+1}$ in terms of ${\rm Vol}_{g,n}$, and for ${\rm
Vol}_{g,0}$ in terms of ${\rm Vol}_{m,k}$ for $m<g$, which, using
our upper bounds on ${\rm Vol}_{g,n}$ for $n>0$, then implies that
${\rm Vol}_{g,0}\le c^gg^{2g}$, and using Penner's lower bounds
for ${\rm Vol}_{g,1}$ also yields ${\rm Vol}_{g,n}\ge A^gg^{2g}$
for all $n$, including the no-puncture case. In particular
Schumacher and Trapani prove (theorem 2 in \cite{schum}) that
\begin{equation}
14\kappa^{3g-3}\ge \kappa^{3g-4}\cdot \sum\limits_{i=0}^{\lfloor
g/2 \rfloor} \de_i. \label{stineq}
\end{equation}

We will now use the Hodge index inequality in the following
form:

\begin{lm}[\cite{demailly}, lemma 5.3, see also \cite{lazars}]
Let $X$ be a complex variety of dimension $N$, and let $E_1\ldots
E_N$ be any nef divisors (i.e. intersecting any effective curve
non-negatively) on $X$. Then for the intersection numbers we have
\begin{equation}
\label{genineq}
\langle E_1\cdots E_N\rangle_X^N\ge\langle
E_1^N\rangle_X\cdots\langle E_N^N\rangle_X
\end{equation}
\end{lm}
If we set $E_1=\ldots=E_p=\alpha$ and $E_{p+1}=\ldots=E_N=\beta$,
this general inequality becomes
\begin{equation}
\langle\alpha^p\cdot\beta^{N-p}\rangle^N\ge\langle\alpha^N\rangle^p
\langle\beta^N\rangle^{N-p}.
\end{equation}

\smallskip
\begin{proof}[Proof of theorem \ref{dbound}]
Applying the last version of Hodge index theorem for $p=1$,
$\alpha=\lambda$ and $\beta=\kappa$ yields
\begin{equation}
\label{appp}
\langle \lambda^{3g-3}\rangle_\omg\langle\kappa^{3g-3}\rangle_\omg^{
3g-4}\le\langle \lambda\kappa^{3g-4}\rangle_\omg^{3g-3}.
\end{equation}
From \cite{wolpert} we recall that in $H^2(\omg)$
\begin{equation}
12\lambda=\kappa+\di+\de_1/2+ \sum_{i>1}\de_i.
\end{equation}
Substituting this in the previous formula, and then using
Schumacher and Trapani's estimate (\ref{stineq}), we get
$$
\langle \lambda\kappa^{3g-4}\rangle_\omg=\frac{1}{12}\left\langle
(\kappa +\di+\de_1/2+\sum\limits_{i>1}\de_i)\kappa^{3g-4}
\right\rangle_\omg<\frac{15}{12}\langle\kappa^{3g-3}\rangle_\omg.
$$
Substituting this in (\ref{appp}) and using the upper bound
of $c^gg^{2g}$ for the Weil-Petersson volume ${\rm
Vol}_{g,0}=\langle\kappa^{3g-3}\rangle/ {(3g-3)!}$, we then
end up with
\begin{equation}
\langle \lambda^{3g-3}\rangle\le\frac{\langle
\lambda\kappa^{3g-4}\rangle^{3g-3}} {\langle
\kappa^{3g-3}\rangle^{3g-4}}<\left(\frac{5}{4}\right)^{3g-3}
\langle \kappa^{3g-3}\rangle<(3g-3)!C^gg^{2g},
\end{equation}
where $C:=5^3c/4^3$ is a new explicit constant.

Combining this with theorem \ref{fformula}, we then obtain the
following bound for the degree
\begin{equation}
\deg Th(\J_g^{2,4})=\int\limits_\mle\Omega^{3g-3}<
(3g-3)!C^gg^{2g}2^{g^2-g+3}\prod\limits_{k=1}^g(2^{2k}-1)
\end{equation}
which is absolutely huge, but explicit and finite nonetheless.
\end{proof}

\section{Eliminating the unknowns in the KP equation}
\label{unknowns} Shiota has obtained a solution to the Schottky
problem via the KP equation for the theta functions:

\begin{thm}[\cite{shiota}]
A principally polarized abelian variety is the Jacobian of a
Riemann surface if and only if there exist three vectors $u\ne
0,v,w\in\C^g$, and a constant $c\in \C$ such that the theta
function $\t(\tau,z):=\tt{0}{0}(\tau, z)$ of this ppav satisfies
the following differential equation (the KP equation) for all
values of $z$:
\begin{equation}
\label{KPShiota}
\t_{uuuu}\t-4\t_{uuu}\t_{u}+3\t_{uu}\t_{uu}+\t_u\t_w-
4\t_{uw}\t+3\t_{vv}\t-3\t_v\t_v+8c\t^2=0,
\end{equation}
where the subscript denotes differentiation with respect to the
$z$ variables in the direction of the vector indicated.
\end{thm}

Using the addition theorem and other properties of theta
functions, this equation was reformulated in terms of theta
constants of the second order by Dubrovin. We will use the
following convention:

\begin{df}
For $u,v\in\C^g$ we denote by
\begin{equation}
uv\p\T[\e](\tau):=\sum\limits_{i,j=1}^gu_iv_j
\frac{\p\T[\e](\tau)} {\p\tau_{ij}}
\end{equation}
the value of the bilinear form determined by the matrix
$\frac{\p}{\p\tau_{ij}}\T[\e]$ at the pair of vectors $(u,v)$. We
also define $uvwx\p^2\T[\e](\tau)$ similarly, as the convolution
of rank four tensors $u\otimes v\otimes w\otimes x$ and
$\frac{\p^2}{\p\tau_{ij}\p\tau_{kl}}\T[\e](\tau)$.
\end{df}

In these notations the KP equation can be reformulated as
follows:

\begin{prop}[\cite{dubrovin}]
\label{dubr} The theta function of a principally polarized abelian
variety with period matrix $\tau$ satisfies the KP equation
(\ref{KPShiota}) if and only if for the same $u,v,w$ and $c$ all
theta constants of the second order satisfy the following
differential equations:
\begin{equation}
\label{KPDubrovin}
u^4\p^2\T[\e](\tau)+(\frac{3}{4}v^2-uw)\p\T[\e](\tau)+c
\T[\e](\tau)=0 \qquad \forall \e.
\end{equation}
\end{prop}

\begin{rem} Dubrovin in \cite{dubrovin} indicates and
Sasaki in \cite{sasaki} proves that the $g(g+1)/2\times 2^g$
matrix $\lbrace\frac{\p}{\p\tau_{ij}}\T[\e], \T[\e]\rbrace$ has
maximal rank, $g(g+1)/2$, for irreducible abelian varieties. It
then follows that if the KP equation in Dubrovin's formulation has
a solution with $u=0$, i.e. if we have $v^2\p\T[\e]+4c\T[\e]/3=0$
valid for all $\e$, then we would need to have $v=0$ and $c=0$.
Thus in the sequel we do not have to worry about the $u\ne 0$
condition in the KP equation.
\end{rem}

The KP equation does not directly allow one to verify whether an
explicitly given principally polarized abelian variety is a
Jacobian. First, one needs to be able to eliminate the unknowns
$u,v,w$ and $c$ from the system of equations (\ref{KPDubrovin}).
Second, even after we eliminate the unknowns, we end up with a
system of non-linear differential equations for theta constants,
the validity of which we need to check at the point $\tau$. We
will now deal with these problems.

To eliminate $u,v,w,$ and $c$ from (\ref{KPDubrovin}), treat
them as unknowns, and all expressions in theta constants and
derivatives --- as given coefficients. Then we have a system
of $M:=2^g$ polynomial equations in $N:=3g+1$ variables:
$\lbrace f_i(x_1,\dots,x_N)=0\rbrace$ for $i=1\ldots M$, with
$\deg f_i=4$, and we ask whether it admits a solution.

The problem can be looked at in two different ways. On one hand,
we can think of this as a problem of eliminating variables
$x_1,\dots,x_N$ from the system of equations. Once we eliminate
all of them, we end up with a system of relations among the
coefficients of $f_i$. On the other hand, instead of doing
elimination we can ask a global question whether this system has a
solution --- this amounts to dealing with Nullstellensatz. The
currently available techniques in Nullstellensatz and elimination
(see \cite{clos}) actually yield similar results, but using
Nullstellensatz is much easier. We use the following version of
the effective Nullstellensatz:

\begin{thm}[\cite{kollar}]
For $f_i\in\C[x_1,\ldots,x_N]$, $\deg f_i\le d$ the system
$\lbrace f_i(x_1,\ldots, x_N)=0\rbrace$, $i=1,\ldots,M$ does not
have a solution in $\C^N$ if and only if there exist
$c_1,\ldots,c_M\in\C[x_1,\ldots, x_N]$ such that $\sum c_if_i=1$
and $\deg (c_if_i)\le d^N\ \ \forall i$.
\end{thm}

Applying this to our situation allows us to effectively
eliminate the variables from the KP equation:

\begin{lm}
Using the effective Nullstellensatz to eliminate the unknowns
in (\ref{KPDubrovin}), we obtain a system of finitely many
parameterless non-linear differential equations and
inequalities for theta constants, which is equivalent to the
KP equation for the theta function.
\end{lm}
\begin{proof}
The total degree of each $f_i$ is equal to 4, and there are
$3g+1$ variables, so by Koll\'ar's result we can choose all
$c_i$ to have total degree at most $4^{3g}$. The number of
monomials in $3g+1$ variables of total degree at most
$4^{3g}$ is $\binom{4^{3g}+3g+1}{3g+1}$. Thus the total
number of variables we have, which is the total number of
undetermined coefficients of all $c_i$ together, is
$K:=2^g\binom{4^{3g}+3g+1}{3g+1}$. The condition $1=\sum
c_if_i$ gives one equation for the coefficient of each
monomial of $\sum c_if_i$ --- there are
$L:=\binom{4^{3g}+3g+5}{3g+1}$ that there can be, as this sum
is a polynomial of total degree at most $4+4^{3g}$. So we get
a system of $L$ linear equations for $K$ unknowns, all but
one of which have no constant term.

Such a system does not have a solution if and only if the one
equation with the non-zero constant term is a linear
combination of all the others. Verifying that this one
equation (call it $R$ --- we think of it as a row of a
matrix) is a linear combination of all the others (call the
matrix of those $A$) amounts to verifying that the rank of
the matrix $A$ together with $R$ is the same as the rank of
$A$ by itself.

Computing the rank by checking the vanishing and non-vanishing of
the determinants of minors of increasing size, this verification
can be done effectively. Thus the condition that the rank of $A$
together with $R$ is the same as the rank of $A$ is just a certain
system of polynomial equations and inequalities for the entries of
$A$ and $R$.
\end{proof}

\section{Solving the Schottky problem effectively}
\label{solving} For the purposes of this section let us fix the
genus/dimension $g$, denote $N:=2^g-1$, and drop $g$ and $2,4$ in
all notations. All polynomials considered are homogeneous.

We are dealing with the following loci: $\J\subset\A$, and
$Th(\J)\subset Th(\A)\subset\P^N$. In the previous sections of
this work we obtained a formula for $L:=\deg Th(\A)$, and an
explicit upper bound $K>\deg Th(\J)$.

We want to determine the algebraic equations defining these
loci in $\P^N$, i.e. the ideals $I(Th(\A)),I(Th(\J))\subset
\C[x_0,\ldots,x_N]^{ \rm hom}$, in the algebra of homogeneous
polynomials in $N+1$ variables. Notice that for any algebraic
subvariety $Z\subset\P^N$ the ideal $I(Z)$ is generated by
its elements of degree at most $\deg Z$. This crude
observation will allow us to effectively obtain the ideals we
are interested in.

We will use the following generalization of Bezout's theorem:
\begin{lm}[\cite{hartshorne}, Theorem 1.7.7]
Let $H$ be a hypersurface in $\P^N$ and let $Z$ be an irreducible
algebraic subvariety of $\P^N$. If $Z$ is not contained in $H$,
let us denote by $\lbrace Z_i\rbrace$ the irreducible components
of $H\cap Z$. Then the following equality holds:
$$
\deg H\cdot \deg Z=\sum\limits_i\operatorname{mult}_{Z_i}(H\cap Z)
\cdot\deg Z_i,
$$
where $\operatorname{mult}_{Z_i}(H\cap Z)$ denotes the
multiplicity with which the intersection $H\cap Z$ contains
$Z_i$.
\end{lm}

What follows from this lemma is that if a polynomial (defining
$H$) vanishes at a point of the subvariety $Z$ to a high enough
order along $Z$ (this is the high multiplicity condition), then it
vanishes identically on $Z$.

\begin{thm}
\label{genabvar} A system of generators of $I(Th(\A))$ can be
obtained effectively.
\end{thm}
\begin{proof}

Let $P\in\C[\T[\e]]_{\le L}$ be a homogeneous polynomial in theta
constants of degree at most $L$ with undetermined coefficients. We
would like to know whether $P\in I(Th(\A))$. Fix some specific
$p\in\H_g$ --- for example, the point of the moduli space
corresponding to the Jacobian of the hyperelliptic curve for which
the period matrix is computed in \cite{schindler}. The
$(L^2+1)$-jet of $Th(\A)$ at $Th(p)$ is the linear span of the set
of all partial derivatives of all theta constants $\T[\e]$ with
respect to all $\tau_{ij}$ to order $L^2+1$, evaluated at $p$ ---
thus it can be computed effectively. The conditions for $P$ to
vanish along $Th(\A)$ at the point $Th(p)$ to order $L^2+1$ are
thus effectively a finite system of algebraic equations for the
coefficients of $P$. Using the effective Nullstellensatz (possible
as the number of equations and their degree are explicit functions
of an explicit $L$) then allows us to effectively choose the
generators for the ideal of all $P$ that satisfy these equations
(see \cite{ein},\cite{kollar}), which will thus be a basis for
$I(Th(\A))$.
\end{proof}

Now for the Schottky problem:
\begin{thm}
A system of generators of $I(Th(\J))$ can be obtained effectively.
\end{thm}
\begin{proof}
In the previous section we effectively reformulated the KP
equation as a finite system of explicit parameterless
polynomial equations and inequalities for theta constants and
their first and second derivatives. Let $S$ denote this
resulting system of polynomial equations in theta constants
and their derivatives. We will not need to look at the
inequalities in detail.

Since the KP gives a solution to the Schottky problem, the locus
$\J\subset\A$ locally near $p$ is the locus of solutions of
$\lbrace S=0\rbrace$ in $\A$ --- the inequalities may only serve
to cut away extra components of the solution set of $\lbrace
S=0\rbrace$ at a certain distance from $p$. Recall that the locus
of Jacobians is irreducible.

To apply the same argument as in theorem \ref{genabvar} we need to
compute the $(K^2+1)$-jet of the set $Th(\lbrace S=0\rbrace)$ at
$Th(p)$ effectively. To do it we need to differentiate the
equations of the system $S$ with respect to some $\T[\e]$ a number
of times (no more than $K^2+1$ times), and then evaluate at $p$.
Noticing that
$$
\left.\frac{\p}{\p\T[\e]}\left(\frac{\p\T[\de]}{\p\tau_{ij}}\right)
\right|_p=\sum_{kl}\left(\left.\frac{\p\T[\e]}{\p\tau_{kl}}
\right|_p\right)^{-1}\cdot\left.\frac{\p^2\T[\de]}{\p
\tau_{kl}\p\tau_{ij}}\right|_p,
$$
we see that the $(K^2+1)$-jet of $Th(\lbrace S=0\rbrace)$ near
$Th(p)$ only depends on the values of theta constants and their
derivatives with respect to $\tau$ at $p$; thus this jet can be
computed since we know all the equations of the system $S$
explicitly, and the values of theta constants and their
derivatives at $p$ are also explicit numbers. Thus the same
argument as in the previous theorem applies, and we have an
effective way of solving the Schottky problem algebraically.
\end{proof}

\section*{Acknowledgements}
The author would like to thank Professor Yum-Tong Siu for
suggesting the problem, teaching me effective complex geometry,
and for inspiration and guidance throught the work on this
project. Many of the ideas used in this work are either directly
due to professor Siu's ingenuous insight, or have stemmed from the
author's discussions with professor Siu. I am grateful to
professor Siu for sharing with me his thorough understanding of
algebraic geometry and complex analysis, and making me aware of
the beauty and the intricacy of the subject.

We are grateful to Carel Faber for discussing with us earlier
versions of this work, providing the intersection numbers needed
to compute the degrees, and bringing to our attention and
explaining to us the work of van der Geer \cite{vdgeer}.

I would also like to thank Corrado De Concini, Ron Donagi, Gavril
Farkas, Tom Graber, Joseph Harris and Sean Keel for sharing their
insights in moduli theory.

This material is based upon work partially supported under a
National Science Foundation Graduate Fellowship, and is partly
based upon part of the author's Ph.D. dissertation at Harvard
University under the direction of professor Yum-Tong Siu.

\end{document}